\input amstex
\input amsppt.sty
\magnification\magstep1

\def\ni\noindent
\def\sbs{\subset}

\def\as{\operatorname{asdim}}

\def\diam{\operatorname{diam}}
\def\dim{\operatorname{dim}}

\def\R{\text{\bf R}}

\def\Q{\text{\bf Q}}
\def\Z{\text{\bf Z}}

\def\G{\Gamma}

\def\Vcal{\Cal V}

\def\Ucal{\Cal U}

\def\mnorm#1{\| #1 \|}
\def\grpgen#1{ \langle #1 \rangle}

\hoffset= 0.0in
\voffset= 0.0in
\hsize=32pc
\vsize=38pc
\baselineskip=24pt
\NoBlackBoxes
\topmatter
\author
A. Dranishnikov and J. Smith
\endauthor

\title
Asymptotic dimension of discrete groups
\endtitle
\abstract
We extend Gromov's notion of asymptotic dimension of finitely generated
groups to all discrete groups. In particular, we extend the Hurewicz type
theorem proven in [B-D2] to general groups. Then we use this extension
to prove a formula for the asymptotic dimension
of finitely generated solvable groups in terms of their Hirsch length.
\endabstract

\thanks The first author was partially supported by NSF grants DMS-0305152
\endthanks

\address University of Florida, Department of Mathematics, P.O.~Box~118105,
358 Little Hall, Gainesville, FL 32611-8105, USA
\endaddress


\email  dranish\@math.ufl.edu 
justins\@math.ufl.edu
\endemail

\endtopmatter

\document

\head \S1 Introduction\endhead

The main purpose of this paper is to extend the formula $\as
\Gamma\le h(\Gamma)$ from [B-D2], which gives an upper bound  for the 
asymptotic dimension of
nilpotent finitely generated groups in terms of the Hirsch length,
to solvable groups. The main problem here is that even working
with finitely generated groups one has to consider infinitely
generated ones. This is because a subgroup of a finitely presented
group can be infinitely generated. Originally, asymptotic dimension
was defined only for finitely generated groups [Gr]. It was
actually defined as a quasi-isometry invariant of metric spaces.
Since all word metrics induced by a finite set of generators in a
finitely generated group are quasi-isometric, the asymptotic
dimension gives a group invariant. It turns out that the
asymptotic dimension of a finitely generated group is a coarse
invariant (definition below) [Ro2].

In this paper we notice that every countable group admits a left
invariant proper metric and any two such metrics are coarsely
equivalent.  Thus, one can extend the notion of asymptotic
dimension to all countable groups. We show in \S2 that the
Hurewicz type formula of [B-D2] still holds for countable groups.
As a result we prove the asymptotic dimension formula for
solvable groups (\S3). The main ingredient of \S2 is Theorem 2.1
which states that the asymptotic dimension of a countable group
$G$ is the maximum of the asymptotic dimensions of the finitely generated
subgroups. This rule leads to a definition of the asymptotic
dimension for an arbitrary group (\S2).

\

We recall that a map $f:(X,d)\to (Y,d')$ between metric spaces is
a {\it coarse embedding} if there are tending to infinity functions
$\rho_1,\rho_2:\R_+\to\R_+$ such that
$$
\rho_1(d(x,x'))\le d'(f(x),f(x'))\le\rho_2(d(x,x'))
$$
for all $x,x'\in X$.
A coarse embedding $f:(X,d)\to (Y,d')$ is a {\it coarse equivalence}
if the image $f(X)$ is $R$-dense in $Y$ for some $R>0$.
A metric $d$ on $X$ is called {\it proper} if every closed ball
$B_r(x,d)=\{y\in X\mid d(x,y)\le r\}$ in $X$ is compact. We denote by
$\dot B_r(x,d)=\{y\in X\mid d(x,y)<r\}$
the open ball of radius $r$ centered at $x$. A metric $d$ on a group $G$ is
called {\it left invariant} if $d(x,y)=d(gx,gy)$ for all $g,x,y\in G$.

\proclaim{Proposition 1.1}
Let $d$ and $d'$ be two proper left invariant metrics on a group $G$ which 
induce the same topology.
Then the identity map $id:(G,d)\to (G,d')$ is a coarse equivalence.
\endproclaim
\demo{Proof} Let $e\in G$ be the unit.
We define
$$
\rho_1(t)=\min\{d'(e,z)\mid z\in G \setminus \dot B_t(e,d)\}
$$
and
$$
\rho_2(t)=\max\{d'(e,z)\mid z\in B_t(e,d)\}.
$$
The properness of $d$ and $d'$ implies that $\rho_1$ (and hence $\rho_2$)
tends to infinity. Let $x,y\in G$. We apply the obvious inequalities
$$
\rho_1(d(e,z))\le d'(e,z)\le\rho_2(d(e,z))
$$
with $z=x^{-1}y$ to obtain the required inequalities
$$
\rho_1(d(x,y))\le d'(x,y)\le\rho_2(d(x,y)).
$$
\qed
\enddemo
\proclaim{Corollary 1.2} Let $\Gamma\subset G$ be a 
closed subgroup and suppose that $\Gamma$ and $G$ are supplied with
proper left invariant metrics $d$ and $d'$, respectively, that induce the 
same topology on $\G$. Then the
inclusion $(\Gamma, d)\to (G,d')$ is a coarse embedding.
\endproclaim

Let $\Gamma$ be a group.  We recall that
a map $||\cdot ||:\Gamma \rightarrow [0, \infty) $
is said to be a norm on $\Gamma$ if the following three properties hold:
\roster
\item{} $\|x\|=0$ if and only if $x=e$,
\item{} $\|x^{-1}\|=\|x\|$ for all $x\in \Gamma$, and
\item{} $\|xy\|\leq \|x\|+\|y\|$ for all $x, y\in \Gamma$.
\endroster
Norms on a group $\Gamma$ are in bijective correspondence with left invariant
metrics:
Given a norm $\|\cdot\|$ on $\Gamma$, define a metric
 by $d(x,y)=\|x^{-1}y\|$, and given a left invariant metric $d$ on
$\Gamma$ define $\|x\|=d(e,x)$.  We call a norm {\it proper} if it is proper 
as a map (where $\Gamma$ has the topology induced by the associated metric).  
Note that proper norms correspond to proper metrics.  

DEFINITION.  Let $\Gamma$ be a countable discrete group.  Let $S$ be a
symmetric generating set (possibly infinite), $S=S^{-1}$,
for $\Gamma$.  A {\it weight function} $w:S\rightarrow [0, \infty)$ on $S$
is any positive, proper function such that $w(s^{-1})=w(s)$ for all $s \in S$.
The properness can essentially be viewed as the requirement that
$\lim w = \infty$.

\proclaim{Proposition 1.3}
Every weight function $w$ on a countable group defines a proper norm
$$\|x\|_w=\inf{ \{\sum_{i=1}^{n}w(s_i)| x=s_1 s_2 \cdots s_n , s_i\in S \} }.$$
\endproclaim
\demo{Proof} Obviously, the infimum in the above definition is a
minimum. Then the conditions (1)-(3) are easy to check.\qed
\enddemo

\head \S2 Asymptotic dimension of general groups\endhead

Gromov defined the asymptotic dimension of a metric space $X$ as follows [Gr].

DEFINITION.  We say that a metric space $X$ has asymptotic dimension $\leq n$ if,
for every $d>0$, there is an $R$ and $n+1$ $d-$disjoint, $R-$bounded families
$\Ucal_0, \Ucal_1, \ldots, \Ucal_n$ of subsets of $X$ such that
$\cup_{i=0}^n \Ucal_i$ is a cover of $X$.

We say that a family $\Ucal$ of subsets of $X$ is $R-$bounded if
$\sup \{ \diam{U} | U \in \Ucal \} \leq R$.  Also, $\Ucal$ is
said to be $d-$disjoint
if $d(x,y) > d$ whenever $x\in U$, $y \in V$, $U \in \Ucal$,
$V \in \Ucal$, and $U \neq V$.

It follows from the definition that $\as X\le \as Y$ for 
a subset $X\subset Y$ taken with the restricted metric.

The notion $\as$ is a coarse invariant (see \cite{Ro2}), i.e.
$\as X=\as X'$ for coarsely equivalent metric spaces.
Since for a coarse embedding $f:X\to Y$ the image $f(X)$ 
is coarsely equivalent to $X$, we obtain $\as X\le \as Y$.

In view of Propositions 1.1
and 1.3 the invariant $\as\Gamma$ is well defined for every
countable group $\Gamma$. By Corollary 1.2, $\as\Gamma'\le\as\Gamma$
for any subgroup $\Gamma'\subset\Gamma$.

The following theorem reduces the study of asymptotic dimension of
countable groups to the asymptotic dimension of finitely generated
groups. \proclaim{Theorem 2.1} Let $G$ be a countable group. Then
$$\as G = \sup \as F, $$ where the supremum varies over finitely
generated subgroups $F$ of $G$.
\endproclaim
\demo{Proof} Fix a weight function $w: G \rightarrow [0, \infty)$;
let $\mnorm{\cdot}$ and $d$
denote the induced norm and metric, respectively.  If $\sup \as F =
\infty$, then $\as G = \infty$,
and we are finished.  We will now assume $\sup \as F < \infty$.
Set $m = \sup \as F$.

Let $d>0$ be given.  Set $T = \{ g \mid w(g) < d \}$ and $F =
\grpgen{T}$.  By the definition of
weight function, $T$ is finite and so $F$ is finitely generated.
Thus, $\as F \leq m$.  So there
exist uniformly bounded, $d-$disjoint families $\Ucal_0, \Ucal_1,
\ldots \Ucal_m$ of subsets of $F$
such that $\cup_i \Ucal_i$ is a cover of $F$.  Let $Z$
be a system of representatives for the partition by cosets $G / F$.  For $0\leq i \leq m$, define
$\Vcal_i = \{ zU \mid z \in Z, U \in \Ucal_i \}$.

It is easy to check that $\Vcal_i$ is uniformly bounded for $i=0,1,\ldots m$, and that
$\cup_i \Vcal_i$ is a cover of $G$.  We check now that $\Vcal_i$ is a $d-$disjoint family.
For suppose $zU \neq z'U'$, where $z,z' \in Z$ and $U,U' \in \Ucal_i$.  Let $x \in zU$ and
$y \in z'U'$.  First suppose $z \neq z'$.  Note that $zU \subset zF$ and $z'U' \subset z'F$, so that
$x \in zF$ and $y \in z'F$.  Since $z \neq z'$, $zF \neq z'F$.  But
$xF=zF$ and $yF=zF$.  So $xF \neq yF$, and hence $x^{-1}y \notin F$.  This means that
$x^{-1}y$ cannot be written as a product of the elements of $T$, and so
$d(x,y)=\mnorm{x^{-1}y} \geq d$ by definition of $\mnorm{\cdot}$.  
Now suppose $z=z'$.  So $x=zu$ and $y=zu'$ for
some $u \in U$ and $u' \in U'$.  Since we must have $U \neq U'$,
$U$ and $U'$ are $d-$disjoint; thus, $d(x,y) =d(zu,zu')=d(u,u') \geq d$.  
It follows that the family $\Vcal_i$ is $d-$disjoint.

Since $d>0$ was arbitrary, $\as G \leq m$.  Equality immediately follows. \qed
\enddemo

REMARK 1.  Since every finitely generated subgroup of $\Q$ is cyclic, 
Theorem 2.1 implies in particular that $\as\Q=1$ (see 
[S] for a direct computation).

Theorem 2.1 could lead to the the following generalization.  

DEFINITION. Let $G$ be an arbitrary group. We define
$$
\as G=\sup\{\as F\mid F\subset G \ \text{is finitely generated}\}.
$$
It is immediate that $\as G \leq \as H$ whenever $G \leqslant H$.  

REMARK 2. Define a coarse structure $\Cal{E}$ on $G$ as follows: $E \in \Cal{E}$ if and only if 
$\{ x^{-1}y \mid (x,y) \in E \}$ is a finite set.  It can be shown that the 
definition of asymptotic dimension above coincides with the asymptotic 
dimension of $G$ equipped with the coarse structure $\Cal{E}$.  It is not hard to see 
that all bounded sets are finite. In the case when $G$ is countable, 
this coarse structure is precisely the bounded coarse structure associated 
with a proper, left invariant metric.  Also, $\Cal{E}$ is generated by the 
family $\{ \Delta_g \mid g \in G \}$, where $\Delta_g = \{ (h, hg) \mid h \in G \}$.  
This shows that $\Cal{E}$ is a natural generalization of Example 2.13 from 
[Ro2].  Using the terminology found in this example, $\Delta_g$ is the 
$G-$saturation of $\{ (1,g) \}$, and $\Cal{E}$ is the coarse structure 
generated by the $G-$saturations of the finite subsets of $G \times G$.    
We refer to [Ro2] for the corresponding  definitions.

REMARK 3.  Even for familiar (uncountable) groups, this new definition of 
asymptotic dimension does not agree with the asymptotic dimension obtained when the 
group is equipped with a left invariant, proper metric.  For example, using a 
formula from Theorem 3.2 
below, we can show that $\as \R = \infty$ and $\as \R / \Z = \infty$.   

Though the results stated below hold true for general discrete groups,
the most important case is the case of countable groups.

We recall a notion of the $R$-stabilizer from [B-D1].

DEFINITION. Suppose that a group $\Gamma$ acts on a metric space $X$
by isometries. Let $x_0\in X$. For every $R>0$ we define the $R$-stabilizer of
$x_0$ as follows: 
$$W_R(x_0)=\{g\in\Gamma\mid d(g(x_0),x_0)\le R\}.$$
We extend to general groups the following theorems 
from [B-D2].

\proclaim{Theorem 2, [B-D2] }
Suppose that a finitely generated group $\Gamma$ acts on a geodesic space $X$
by isometries. Let $x_0\in X$ and suppose $\as W_R(x_0)\le n$ for all
$R>0$. Then $\as \Gamma\le n+\as X$.
\endproclaim

\proclaim{Theorem 7, [B-D2] (Hurewicz Type Formula)}
Let $\phi:G\to H$ be a homomorphism of of finitely generated
groups with kernel $K$. Then
$$
\as G\le\as H+\as K.
$$
\endproclaim

Since the notion of asymptotic dimension at time when [D-B2] was wirtten
was defined only for finitely generated groups, $K$ is treated in Theorem 7
as metric space
with the word metric restricted from $G$. 

\proclaim{Theorem 2.2 }
Suppose that a group $\Gamma$ acts on a geodesic space $X$
by isometries. Let $x_0\in X$ and suppose $\as W_R(x_0)\le n$ for all
$R>0$. Then $\as \Gamma\le n+\as X$.
\endproclaim
\demo{Proof}
Let $F\subset \Gamma$ be a finitely generated subgroup. Then the $R$-stabilizer
with respect to $F$ is contained in the $R$-stabilizer
with respect to $\Gamma$. Then by Theorem 2 of [B-D2],
$\as F\le n+\as X$. By the definition, $\as \Gamma\le n+\as X$.\qed
\enddemo
\proclaim{Theorem 2.3 (Hurewicz Type Formula)}
Let $\phi:G\to H$ be a homomorphism of groups with kernel $K$. Then
$$
\as G\le\as H+\as K.
$$
\endproclaim
\demo{Proof}
Let $F\subset G$ be a finitely generated subgroup and
let $d^F$ denote the word metric on $F$.
Then $$\as F\le \as \phi (F)+\as (K\cap F,d^F|_{K\cap F})$$ by 
Theorem 7 of [B-D2]. If $\as (K\cap F,d^F|_{K\cap F})<\infty$, then
by Theorem 2.1, $\as (K\cap F,d^F|_{K\cap F})=\as F'$
for some finitely generated group $F'\subset K\cap F$. Then by the definition
$\as \phi (F)\le\as H$ and $\as F'\le\as K$. Hence $\as F \le \as H + \as K$.  
The assertion now follows from the definition of $\as G$.
If $\as (K\cap F,d^F|_{K\cap F})=\infty$, by Theorem 2.1 then there is 
a sequence $F'_n\subset K\cap F$ of finitely generated subgroups with
$\as F'_n\ge n$. Then $\as K=\infty$ and the inequality follows.
\qed
\enddemo

\head \S3  Asymptotic dimension of solvable groups\endhead

A group $T$ is called {\it torsion group} if every element of $T$ has
a finite order.

In view of the facts that every finite group has asymptotic
dimension 0 and every finitely generated
torsion abelian group is finite, we obtain the following.  
\proclaim{Lemma 3.1}
$\as T=0$ for every torsion abelian group.
\endproclaim

\proclaim{Theorem 3.2} Let $A$ be an abelian group. Then
$\as A=\dim_{\Q}(A\otimes\Q).$
\endproclaim
\demo{Proof}
Let $\dim_{\Q}(A\otimes\Q)=n$.
By induction on $n$ we show that $\Z^n\subset A$.
Consider a short exact sequence
$$
0\to\Z\to A @>\psi>> A/\Z\to 0
$$
and tensor multiply it by $\Q$ to obtain that
$\dim_{\Q}((A/\Z)\otimes\Q)=n-1$. By the induction assumption
$\Z^{n-1}\subset A/\Z$. Then a section of $\psi$ over $\Z^{n-1}$ together with
the imbedding of the fiber $Z\to A$ defines a monomorphism
$Z^n\to A$.
Therefore, $\as A\ge n$.

On the other hand since every torsion free $\Z$-module is flat,
we obtain
$$
A/TorA=(A/TorA)\otimes\Z\subset (A/TorA)\otimes\Q= A\otimes\Q=\Q^n.
$$
By Remark 1 and the product theorem for asymptotic dimension [D-J]
we have $\as\Q^n=n$.
Hence $\as (A/Tor A)\le n$. Then by the Hurewicz type formula (Theorem 2.3)
and Lemma 3.1,
$\as A\le \as(A/Tor A)+\as(TorA)\le n+0$.
\qed
\enddemo

\proclaim{Corollary 3.3}  For any
short exact sequence of abelian groups
$$ 0 \to B \to A \to C \to 0 $$
we have that $\as A = \as B + \as C$.  In particular,
$\as A \times B = \as A + \as B$.
\endproclaim

DEFINITION. Let $\Gamma$ be solvable and let
$1=G_0\subset G_1\subset\dots\subset G_n=\Gamma$ be
the commutator series, $G_i=[G_{i+1}, G_{i+1}]$. Then the
{\it  Hirsch length}
of $\Gamma$ is defined as
$$h(\Gamma)=\sum\dim_{\Q}((G_{i+1}/G_i)\otimes\Q).$$

\proclaim{Theorem 3.4}
For a solvable group $\Gamma$, $\as\Gamma\le h(\Gamma)$.
\endproclaim
\demo{Proof}
Apply the Hurewicz type formula and Theorem 3.2\qed
\enddemo

For a virtually polycyclic group $\G$, we define the Hirsch length 
$h(\G)$ to be $h(\G')$, where $\G'$ is a polycyclic, 
normal subgroup of $\G$ of finite index.  It is easy to 
check that this Hirsch length is well defined.   

\proclaim{Theorem 3.5}
For a virtually polycyclic group $\Gamma$, 
$$\as\Gamma = h(\Gamma).$$
\endproclaim
\demo{Proof}
Since both the asymptotic dimension and the Hirsch length are left 
unchanged by passing from a virtually polycyclic group to a 
polycyclic subgroup of finite index, it suffices to prove this 
theorem when $\G$ is polycyclic.  In view of Theorem 3.4 it suffices 
to show the inequality
$\as\Gamma\ge h(\Gamma)$.  By
Theorems 3.1 and  4.28 in [Ra], there is a normal subgroup $H\subset\Gamma$ of
finite index which is a cocompact lattice in a solvable simply
connected Lie group $G$. Then $H$ is  coarsely isomorphic to $G$
taken with an invariant metric. On the other hand the homogeneous
space $G/K$ is coarsely isomorphic to $G$ where $K$ is a maximal
compact subgroup in $G$. The metric on $G/K$ can be taken to be
the Hausdorff metric on the orbits $Kx$. It is well known that the
homogeneous space $G/K$ is homeomorphic to $\R^n$. 
Since every locally compact contractible (metrically) homogeneous metric space
is uniformly contractible, the space $G/K$ is uniformly contractible. 
One can extract from the
proof of Theorem 4.28 in [Ra] that $n=h(\Gamma)$. By a theorem of
Roe, $HX^n(G/K)=H^n_c(G/K)=\Z$ since $G / K$ is a uniformly contractible
metric space, where $HX^*$ denotes the coarse cohomology [Ro1].
Since the coarse cohomology is a coarse invariant,
$HX^n(\Gamma)=\Z$. We show that $\as\Gamma\ge n$.

Assume that $\as\Gamma\le n-1$. Then $\Gamma$ admits an
anti-\v Cech approximation by $n-1$-dimensional locally finite
polyhedra [Dr],[Ro1].
It means there is a direct system of $n-1$-dimensional locally finite 
polyhedra $\{K_i,\phi^i_{i+1}\}$
associated with $\Gamma$
such that
$$
0\to {\lim_{\leftarrow}}^1H^{n-1}_c(K_i)\to HX^n(\Gamma)\to\lim_{\leftarrow}H^n_c(K_i)\to 0.
$$
Being isomorphic to the \v Cech cohomology group of a compact
space every group $H^{n-1}_c(K_i)$ is countable.
Since the lim-one of countable groups is either 0 or uncountable
[Ha] we obtain a contradiction: $\Z=\lim^1_{\leftarrow}H^{n-1}_c(K_i)$.
 \qed
\enddemo

\

\enddocument